\documentclass[12pt]{amsart}

\usepackage{a4,amsmath,amssymb}

\usepackage{tikz}
\usetikzlibrary{arrows}

\usepackage{enumitem}

\setlist{itemsep=1mm,topsep=1mm}
\setitemize{leftmargin=7mm,labelsep=2mm}

\allowdisplaybreaks

\theoremstyle{plain}
\newtheorem*{thmmain}{Main Theorem}
\newtheorem{lemma}{Lemma}

\DeclareMathOperator{\N}{NF}

\DeclareMathOperator{\Cay}{\Gamma}
\DeclareMathOperator{\Cayu}{\Gamma^\prime }

\newcommand{\nf}[1]{\textsf{#1}}

\date{\today}

\begin{document}

\title{Finite presentability and isomorphism of Cayley graphs of monoids}

\author{J. Awang}
\address{School of Mathematics and Statistics, University of St Andrews, St Andrews, Scotland, U.K.}
\email{jsa23@st-andrews.ac.uk}

\author{M. Pfeiffer}
\address{School of Computer Science, University of St Andrews, St Andrews, Scotland, U.K.}
\email{markus.pfeiffer@st-andrews.ac.uk}

\author{N. Ru\v{s}kuc}
\address{School of Mathematics and Statistics, University of St Andrews, St Andrews, Scotland, U.K.}
\email{nik.ruskuc@st-andrews.ac.uk}

\begin{abstract}
Two finitely generated monoids are constructed, one finitely presented the other not, whose (directed, unlabelled) Cayley graphs are isomorphic.
\end{abstract}


\maketitle

\section{Introduction}
\label{sec1}

Is the property of being finitely presented a quasi-isometry invariant for the class of all monoids?
This question is posed by  Gray and Kambites in \cite[Question 1]{graykam13a}.
It is well known that the answer is positive for groups \cite[Proposition I.8.24]{bridsonhaef},
and Gray and Kambites also establish it for monoids with finitely many left- and right ideals 
\cite[Theorem 4]{graykam13a}, and for left cancellative monoids \cite[Theorem A]{graykam13b}.
They also show that finite presentability is not invariant when the direction of edges in Cayley graphs is suppressed (\cite[Theorem 11]{graykam13b}).
The purpose of this paper is to answer the general question negatively, by means of a concrete example.
In fact, more strikingly,
the two monoids $M$ and $N$ that are constructed, one finitely presented and the other not, have 
\emph{isomorphic}
(directed, unlabelled) Cayley graphs.

For a monoid $M$ with a generating set $A$ the (right) \emph{Cayley graph}
$\Cay(M,A)$ is a labelled directed graph with vertices $M$ and, for every $a\in A$, edges
$(x,xa)$ ($x\in M$) labelled $a$.
It encodes the right regular representation of $M$, and hence completely determines $M$
by Cayley's Theorem (for monoids).
One can, of course, dually define the \emph{left} Cayley graph of $M$, which,
unlike for groups, need not be isomorphic to the right one.
Cayley graphs are fundamental objects in combinatorial theory of monoids and groups.
Considering undirected paths in $\Cay(M,A)$ turns $M$ into a metric space; this transition has an enormous significance in group theory, and is essentially the interface between combinatorial and geometric theories of groups; see for instance the Introduction to \cite{delaharpe00}.
In a series of papers \cite{graykam11,graykam13a,graykam13b} Kambites and Gray argue that the `correct' analogue in the case of monoids are semimetric spaces, obtained by considering \emph{directed} paths in the Cayley graph.
In either case, quasi-isometry is defined as a natural geometrical notion, which identifies spaces which `look alike' from `far away'; see \cite[Definition 8.14]{bridsonhaef} and \cite[Definition 7]{graykam13a}.
 The question then arises which algebraic properties of the original monoid $M$ remain `encoded' in the new geometric setting, or, more precisely, are invariant under quasi-isometries.

For the purposes of the present paper we will not actually require any of the geometric apparatus, only the notion of the \emph{unlabelled} Cayley graph
$\Cayu(M,A)$, obtained from $\Cay(M,A)$ by simply `omitting'  all the edge labels. Note that this allows for loops and multiple edges between pairs of vertices, although they will actually not occur in the examples we are about to exhibit.
Our main result is as follows:

\begin{thmmain}
Let $M$ and $N$ be the monoids defined by the presentations
\begin{align*}
&M=\langle a,b \:|\: ab^na=aba\ (n=2,3,4,\dots)\rangle\\
&N=\langle c,d\:|\: cdc=cd^2c=cd^4=cd^3c^2=cd^3cdc\rangle .
\end{align*}
Then
\begin{enumerate}[label=\textup{(\alph*)}, leftmargin=10mm,labelsep=2mm]
\item
$M$ is not finitely presented.
\item
$N$ is finitely presented.
\item
The unlabelled (directed, right) Cayley graphs $\Cayu(M,\{a,b\})$ and $\Cayu(N,\{c,d\})$
are isomorphic.
\end{enumerate}
\end{thmmain}

It is obvious that $N$ is finitely presented, and it is easy to prove that $M$ is not.
Indeed, if it were, then it would be defined by a finite subset of the given presentation,
say
\[
M=\langle a,b \:|\: ab^2a=aba,ab^3a=aba,\dots, ab^{n_0}a=aba\rangle;
\]
this general algebraic fact is simply a way of saying that congruence generation is an algebraic closure operator
\cite[Theorem 5.5]{Burris} and
 is stated explicitly for semigroups in \cite[Theorem 9.14]{CliffordPreston}.
But then no defining relation could be applied to the word $ab^{n_0+1}a$, contradicting the fact that $ab^{n_0+1}a=aba$ in $M$.

Thus, the only task is to prove (c).
We do this by utilising rewriting systems to establish a set of normal forms for each monoid
(Section \ref{sec2}), and then analysing the action of generators on these normal forms.
We first present an intuitive graphical descriptions of the Cayley graphs of $M$ and $N$
(Section \ref{sec3}), and then write down an explicit isomorphism of unlabelled versions
(Section \ref{sec4}).
Since the edge labels play no role in determining the (semi)metric structure on a monoid
(see \cite[Section 2]{graykam13a}), the negative answer to \cite[Question 1]{graykam13a} is an immediate consequence of our theorem.

\section{Normal forms}
\label{sec2}

In order to write down normal forms for $M$ and $N$ we will use rewriting systems, and we begin with a brief overview; for full details see \cite{book93}.
A \emph{rewriting system} $R$ over an alphabet $A$ is a subset of $A^\ast\times A^\ast$, where $A^\ast$ denotes the free monoid over $A$, including the empty word $\epsilon$.
The \emph{single step rewriting relation} $\rightarrow$ on $A^\ast$ is 
$\{ (xuy,xvy)\::\: (u,v)\in R,\ x,y\in A^\ast\}$, and the \emph{rewriting relation}
$\rightarrow^\ast$ is the reflexive, transitive closure of $\rightarrow$.
The equivalence relation generated by $\rightarrow^\ast$ is denoted by $\leftrightarrow^\ast$;
this is precisely the congruence on $A^\ast$ arising from the presentation
$\langle A\:|\: R\rangle$.

The rewriting system $R$ is \emph{noetherian} if there are no infinite sequences of single step rewritings
$w_1\rightarrow w_2\rightarrow w_3\rightarrow \dots$.
A word $w\in A^\ast$ is \emph{irreducible} if there is no word $u$ with $w\rightarrow u$.
If $R$ is noetherian, for every $w\in A^\ast$ there exists an irreducible word $u$
such that $w\rightarrow^\ast u$.

The rewriting system $R$ is \emph{confluent} (respectively locally confluent) if for any $w,u,v\in A^\ast$ with
$w\rightarrow^\ast u$, $w\rightarrow^\ast v$ 
(respectively $w\rightarrow u$, $w\rightarrow v$) there exists $z\in A^\ast$ such that
$u\rightarrow^\ast z$, $v\rightarrow^\ast z$.
For a noetherian rewriting system confluence is equivalent to local confluence \cite[Theorem 1.1.13]{book93}. In turn they are equivalent to the \emph{resolution of overlaps}, namely
for for any $(u,v),(z,t)\in R$ the following being true:
\begin{itemize}
\item[1.]
If $u=pq$, $z=qr$ with $q\neq \epsilon$ then there exists $s\in A^\ast$ such that
$vr\rightarrow^\ast s$ and $pt\rightarrow^\ast s$.
\item[2.]
If $u=pzq$ then there exists $s\in A^\ast$ such that
$v\rightarrow^\ast s$ and $ptq\rightarrow^\ast s$;
\end{itemize}
see \cite[Section 2.3]{book93}.
Situations described under 1) and 2) can be visualised as follows:
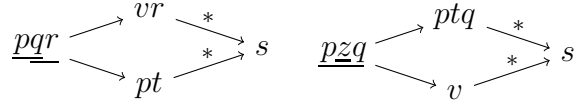
\[
\begin{tikzpicture}[node/.style={text height=5pt,text depth=1pt}]
\node (pqr)  at ( 0,1) {$pqr$\hspace{-6.3mm}\rule[-4pt]{4mm}{0.4pt} \hspace{-3mm}\rule[-5.5pt]{3.8mm}{0.4pt} };
\node (vr)  at ( 1.5,1.5)  {$vr$};
\node (pt)  at ( 1.5,0.5)  {$pt$};
\node (s)  at ( 3,1) {$s$};
\draw[->] (pqr) -- (vr);
\draw[->] (pqr) -- (pt);
\draw[->] (vr) -- (s) node [pos=0.15, above=2pt, right=2pt] {\scriptsize $\ast$};
\draw[->] (pt) -- (s) node [pos=0.8, above=2pt, left=2pt] {\scriptsize $\ast$};
\end{tikzpicture}
\ \ \ 
\begin{tikzpicture}[node/.style={text height=5pt,text depth=1pt}]
\node (pzq)  at ( 0,1) {$pzq$\hspace{-4.2mm}\rule[-2pt]{2mm}{0.4pt} \hspace{-5.7mm}\rule[-4pt]{6mm}{0.4pt}};
\node (ptq)  at ( 1.5,1.5)  {$ptq$};
\node (v)  at ( 1.5,0.5)  {$v$};
\node (s)  at ( 3,1) {$s$};
\draw[->] (pqr) -- (ptq);
\draw[->] (pqr) -- (v);
\draw[->] (ptq) -- (s) node [pos=0.15, above=2pt, right=2pt] {\scriptsize $\ast$}; 
\draw[->] (v) -- (s) node [pos=0.8, above=2pt, left=2pt] {\scriptsize $\ast$};
\end{tikzpicture}
\]
A noetherian and confluent rewriting system is said to be \emph{complete}.
For such a system every equivalence class of $\leftrightarrow^\ast$ contains a unique irreducible element.

It is now easy to check that the presentations defining $M$ and $N$, if interpreted as rewriting systems
\begin{align*}
R_M &: \ ab^na\rightarrow aba\ (n=2,3,\dots)\\
R_N &:\ cd^2c\rightarrow cdc,\ cd^4\rightarrow cdc,\ cd^3c^2\rightarrow cdc,\ cd^3cdc\rightarrow cdc,
\end{align*}
are complete.
Indeed, 
they are noetherian as all rewriting rules are length-reducing.
Furthermore, every overlap in either of the two rewriting systems has the form \mbox{$uxv$\hspace{-6.3mm}\rule[-2pt]{4mm}{0.4pt}\hspace{-2mm}\rule[-3.5pt]{4mm}{0.4pt}}\hspace{1pt}, where $x\in\{a,c\}$ and $ux\rightarrow xyx$ and $xv\rightarrow xyx$ are rewriting rules for an appropriate $y\in\{b,d\}$, and it resolves as follows:
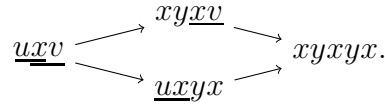
\[
\begin{tikzpicture}[node/.style={text height=5pt,text depth=1pt}]
\node (l)  at ( 0,1) {$uxv$\hspace{-6.5mm}\rule[-2pt]{4mm}{0.4pt}\hspace{-2mm}\rule[-3.5pt]{4.5mm}{0.4pt}};
\node (u)  at ( 2,1.5)  {$xy\underline{xv}$};
\node (d)  at ( 2,0.5)  {$\underline{ux}yx$};
\node (r)  at ( 4,1) {$xyxyx$.};
\draw[->] (l) -- (u);
\draw[->] (l) -- (d);
\draw[->] (u) -- (r);
\draw[->] (d) -- (r);
\end{tikzpicture}
\]

It follows that sets of normal forms are provided by the irreducible words with respect to the two rewriting systems.
In the case of $R_M$ these are precisely the words which do not contain subwords $ab^na$ for $n\geq 2$.
For future use we will record them as follows:

\begin{lemma}
\label{lemma1}
Let
\[
U_M=\{ a^{i_0}ba^{i_1}b\dots ba^{i_{k-1}} ba^{i_k}\::\: k\geq 0,\ i_0,\dots,i_k\geq 1\}.
\]
The monoid $M$ admits the following set of normal forms
\[
\N_M=\{ b^s,b^su,b^sub^t\::\: u\in  U_M,\ s\geq 0,\ t>0\}.
\]
\end{lemma}

The irreducible normal forms for $N$ are just a little bit more complicated.
Clearly, all $c^m$, $d^m$ ($m\geq 0$) are irreducible.
Let us consider an arbitrary irreducible word $w$ which contains occurrences of both $c$ and $d$,
and write $w$ in the following form:
\begin{gather*}
w=d^{j_0} c^{i_1} d^{j_1}\dots c^{i_{k-1}} d^{j_{k-1}}c^{i_k} d^{j_k}\\
(k\geq 1;\ i_1,\dots,i_k\geq 1;\ j_0,j_k\geq 0;\ j_1,\dots,j_{k-1}\geq 1).
\end{gather*}
Bearing in mind the rewriting rules from $R_N$ we see that the following hold:
\begin{itemize}
\item
None of $j_1,\dots, j_k$ can exceed $3$ (because of $cd^4\rightarrow cdc$).
\item
None of $j_1,\dots,j_{k-1}$ equal $2$ (because of $cd^2c\rightarrow cdc$).
\item
If $j_l=3$ for some $l=1,\dots,k-1$, then for the subsequent indices we must have
$
i_{l+1}=\dots=i_k=1,\ j_{l+1}=\dots=j_{k-1}=3
$
(because of $cd^3c^2\rightarrow cdc$ and $cd^3cdc\rightarrow cdcdc$).
\end{itemize}
We can summarise our findings as follows:

\begin{lemma}
\label{lemma2}
Let
\[
 U_N=\{ c^{i_0}dc^{i_1}d\dots dc^{i_{k-1}} dc^{i_k}\::\: k\geq 0,\ i_0,\dots,i_k\geq 1\}.
\]
The monoid $N$ admits the following set of normal forms
\[
\N_N=\{ d^p,d^pu,d^pu(d^3c)^qd^r\::\: u\in  U_N,\ p,q\geq 0,\ 0\leq r\leq 3,\ q+r>0\}.
\]
\end{lemma}

\section{Construction of Cayley graphs}
\label{sec3}

In this section we describe what the Cayley graphs for $M$ and $N$ look like.
Each construction will be in three steps.
Since these steps turn out to be identical, as are the starting ingredients, provided we ignore the edge labels,
this provides an intuitive explanation of the isomorphism between $\Cayu(M,\{a,b\})$ and $\Cayu(N,\{c,d\})$.
The actual formal proof will follow in Section \ref{sec4}.

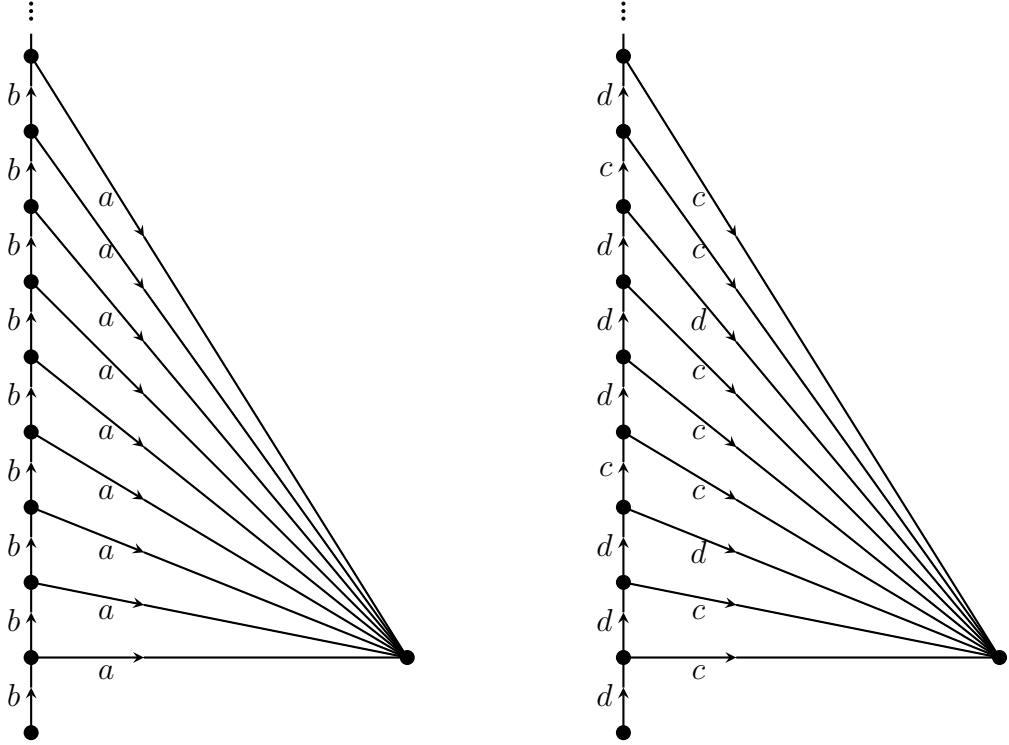
\begin{figure}

\begin{center}

\begin{tikzpicture}


\begin{scope}[radius=1mm]
\foreach \x in {0,...,9}
{\fill (0,\x) circle;}
\fill (5,1) circle ;
\end{scope}

\foreach \i in {1,...,9} 
{\draw (0,\i) -- (5,1);}

\node at (1,0.8) {$a$};
\node at (1,1.6) {$a$};
\node at (1,2.4) {$a$};
\node at (1,3.2) {$a$};
\node at (1,4) {$a$};
\node at (1,4.8) {$a$};
\node at (1,5.5) {$a$};
\node at (1, 6.4) {$a$};
\node at (1, 7.1) {$a$};

\foreach \i in {0.5,1.5,...,8.5}
{\node at (0,\i) [left] {$b$};}

\begin{scope}[thick]
\foreach \i in {1,...,9}
{\draw[-stealth] (0,\i - 1) -- (0,\i - 0.4);
\draw (0,\i - 0.4) -- (0,\i);
\draw[-stealth] (0,\i,0) -- (1.5,\i - 1.5*\i /5 + 1.5/5);
\draw (1.5,\i - 1.5*\i /5 + 1.5/5) -- (5,1);
}
\draw (0,9)--(0,9.3);
\end{scope}

\draw[fill] (0,9.5) circle [radius=0.2mm];
\draw[fill] (0,9.6) circle [radius=0.2mm];
\draw[fill] (0,9.7) circle [radius=0.2mm];

\end{tikzpicture}
\hspace{20mm}
\begin{tikzpicture}


\begin{scope}[radius=1mm]
\foreach \x in {0,...,9}
{\fill (0,\x) circle;}
\fill (5,1) circle ;
\end{scope}

\foreach \i in {1,...,9} 
{\draw (0,\i) -- (5,1);}

\node at (1,0.8) {$c$};
\node at (1,1.6) {$c$};
\node at (1,2.4) {$d$};
\node at (1,3.2) {$c$};
\node at (1,4) {$c$};
\node at (1,4.8) {$c$};
\node at (1,5.5) {$d$};
\node at (1, 6.4) {$c$};
\node at (1, 7.1) {$c$};

\foreach \i in {0.5,1.5,2.5,4.5,5.5,6.5,8.5}
{\node at (0,\i) [left] {$d$};}

\node at (0,3.5) [left] {$c$};
\node at (0,7.5) [left] {$c$};

\begin{scope}[thick]
\foreach \i in {1,...,9}
{\draw[-stealth] (0,\i - 1) -- (0,\i - 0.4);
\draw (0,\i - 0.4) -- (0,\i);
\draw[-stealth] (0,\i,0) -- (1.5,\i - 1.5*\i /5 + 1.5/5);
\draw (1.5,\i - 1.5*\i /5 + 1.5/5) -- (5,1);
}
\draw (0,9)--(0,9.3);
\end{scope}

\draw[fill] (0,9.5) circle [radius=0.2mm];
\draw[fill] (0,9.6) circle [radius=0.2mm];
\draw[fill] (0,9.7) circle [radius=0.2mm];

\end{tikzpicture}

\caption{The graphs $\Gamma_M^{(1)}$ (on the left) and $\Gamma_N^{(1)}$ (on the right).}
\label{fig1}
\label{fig4}
\end{center}
\end{figure}

\begin{figure}

\begin{center}

\begin{tikzpicture}[x=1cm, y=1cm]


\begin{scope}[thick]

\draw [fill=lightgray] (1,1) -- (2,1) -- (1,5) -- (1,1);
   \draw [fill=lightgray] (2,2) -- (3,2) -- (2,5) -- (2,2);
       \draw [fill=lightgray] (3,3) -- (4,3) -- (3,5) -- (3,3);
   \draw [fill=lightgray] (5,2) -- (6,2) -- (5,5) -- (5,2);
\draw [fill=lightgray] (7,1) -- (8,1) -- (7,5) -- (7,1);
   \draw [fill=lightgray] (8,2) -- (9,2) -- (8,5) -- (8,2);
       \draw [fill=lightgray] (9,3) -- (10,3) -- (9,5) -- (9,3);
   \draw [fill=lightgray] (11,2) -- (12,2) -- (11,5) -- (11,2);
\draw [fill=lightgray] (13,1) -- (14,1) -- (13,5) -- (13,1);
\draw [fill=lightgray] (13,1) -- (14,1) -- (13,5) -- (13,1);

\draw[-stealth] (0,0)--(0.6,0);
\draw (0.6,0)--(1,0);
\draw[-stealth] (1,0)--(4,0);
\draw (4,0)--(7,0);
\draw[-stealth] (7,0)--(10,0);
\draw (10,0)--(13,0);
\draw[-stealth] (1,1)--(1.6,1);
\draw (1.6,1)--(2,1);
\draw[-stealth] (2,1)--(3.5,1);
\draw (3.5,1)--(5,1);
\draw[-stealth] (7,1)--(7.6,1);
\draw (7.6,1)--(8,1);
\draw[-stealth] (8,1)--(9.5,1);
\draw (9.5,1)--(11,1);
\draw[-stealth] (13,1)--(13.6,1);
\draw (13.6,1)--(14,1);
\draw[-stealth] (2,2)--(2.6,2);
\draw (2.6,2)--(3,2);
\draw[-stealth] (5,2)--(5.6,2);
\draw (5.6,2)--(6,2);
\draw[-stealth] (8,2)--(8.6,2);
\draw (8.6,2)--(9,2);
\draw[-stealth] (11,2)--(11.6,2);
\draw (11.6,2)--(12,2);
\draw[-stealth] (3,3)--(3.6,3);
\draw (3.6,3)--(4,3);
\draw[-stealth] (9,3)--(9.6,3);
\draw (9.6,3)--(10,3);

\draw[-stealth] (1,0)--(1,0.6);
\draw (1,0.6)--(1,1);
\draw[-stealth] (7,0)--(7,0.6);
\draw (7,0.6)--(7,1);
\draw[-stealth] (13,0)--(13,0.6);
\draw (13,0.6)--(13,1);
\draw[-stealth] (2,1)--(2,1.6);
\draw (2,1.6)--(2,2);
\draw[-stealth] (5,1)--(5,1.6);
\draw (5,1.6)--(5,2);
\draw[-stealth] (8,1)--(8,1.6);
\draw (8,1.6)--(8,2);
\draw[-stealth] (11,1)--(11,1.6);
\draw (11,1.6)--(11,2);
\draw[-stealth] (3,2)--(3,2.6);
\draw (3,2.6)--(3,3);
\draw[-stealth] (9,2)--(9,2.6);
\draw (9,2.6)--(9,3);

\newcommand{\straightend}[2]{
\draw (#1,#2) -- (#1+1.4,#2);
\draw (#1+1.5,#2) circle [radius=0.1mm];
\draw (#1+1.6,#2) circle [radius=0.1mm];
\draw (#1+1.7,#2) circle [radius=0.1mm];
}
\straightend{13}{0}
\straightend{11}{1}
\straightend{9}{2}
\straightend{5}{1}
\straightend{3}{2}

\newcommand{\cornerend}[2]{
\draw (#1,#2) -- (#1 + 0.3, #2);
\draw (#1,#2) -- (#1 , #2 +0.4); (14,1) -- (14,1.4);
\draw (#1 + 0.2,#2 + 0.2) circle [radius=0.1mm];
\draw (#1 + 0.3,#2 + 0.3) circle [radius=0.1mm];
\draw (#1 + 0.4,#2 + 0.4) circle [radius=0.1mm];
}
\cornerend{14}{1}
\cornerend{12}{2}
\cornerend{12}{2}
\cornerend{12}{2}
\cornerend{10}{3}
\cornerend{6}{2}
\cornerend{4}{3}

\end{scope}

\begin{scope}[radius=1mm]

\fill (0,0) circle;
\fill (1,0) circle;
\fill (7,0) circle;
\fill (13,0) circle;
\fill (1,1) circle;
\fill (2,1) circle;
\fill (5,1) circle;
\fill (7,1) circle;
\fill (8,1) circle;
\fill (11,1) circle;
\fill (13,1) circle;
\fill (14,1) circle;
\fill (2,2) circle;
\fill (3,2) circle;
\fill (5,2) circle;
\fill (6,2) circle;
\fill (8,2) circle;
\fill (9,2) circle;
\fill (11,2) circle;
\fill (12,2) circle;
\fill (3,3) circle;
\fill (4,3) circle;
\fill (9,3) circle;
\fill (10,3) circle;
\end{scope}

\begin{scope}[above=-2pt]
\node at (0.5,0)  {$a,c$};
\node at (4,0)  {$a,c$};
\node at (10,0)  {$a,c$};
\node at (1.5,1) [above] {$a,c$};
\node at (3.5,1) [above] {$a,c$};
\node at (7.5,1) [above] {$a,c$};
\node at (9.5,1) [above] {$a,c$};
\node at (13.5,1) [above] {$a,c$};
\node at (2.5,2) [above] {$a,c$};
\node at (5.5,2) [above] {$a,c$};
\node at (8.5,2) [above] {$a,c$};
\node at (11.5,2) [above] {$a,c$};
\node at (3.5,3) [above] {$a,c$};
\node at (9.5,3) [above] {$a,c$};
\end{scope}

\node at (1,0.5) [right] {$b,d$};
\node at (7,0.5) [right] {$b,d$};
\node at (13,0.5) [right] {$b,d$};
\node at (2,1.5) [right] {$b,d$};
\node at (5,1.5) [right] {$b,d$};
\node at (8,1.5) [right] {$b,d$};
\node at (11,1.5) [right] {$b,d$};
\node at (3,2.5) [right] {$b,d$};
\node at (9,2.5) [right] {$b,d$};

\end{tikzpicture}

\caption{The graphs $\Gamma_M^{(2)}$ (reading the $a$, $b$ labels) and $\Gamma_N^{(2)}$ (reading the $c$, $d$ labels). Each constituent copy of $\Gamma_M^{(1)}$ (resp. $\Gamma_N^{(1)}$) is indicated by its initial $ba$ (resp. $dc$) path and a shaded triangle.}
\label{fig2}
\label{fig5}
\end{center}
\end{figure}
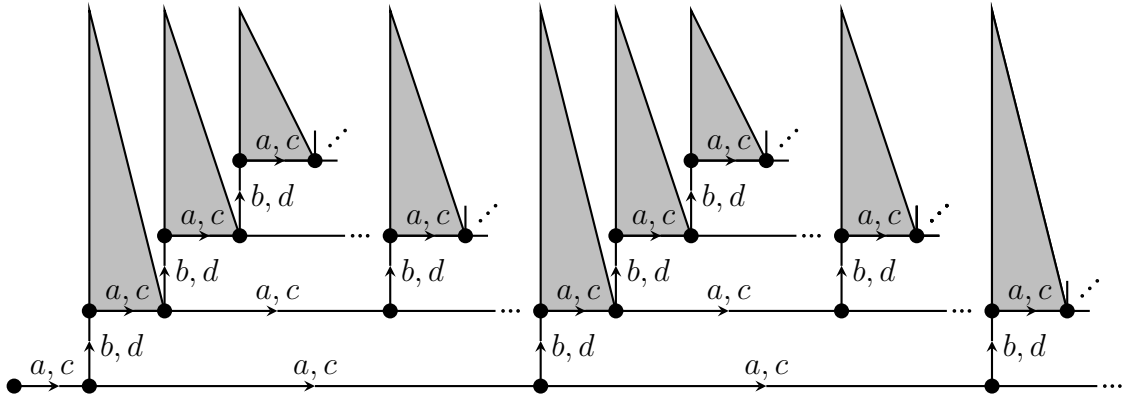

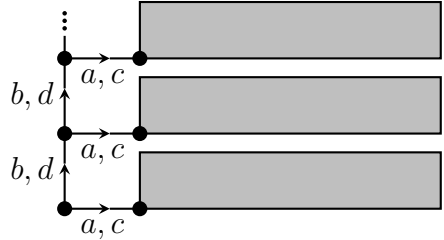
\begin{figure}

\begin{center}

\begin{tikzpicture}


\begin{scope}[thick]

\draw [fill=lightgray] (1,0) -- (5,0) -- (5,0.75) -- (1,0.75)--(1,0);
\draw [fill=lightgray] (1,1) -- (5,1) -- (5,1.75) -- (1,1.75)--(1,1);
\draw [fill=lightgray] (1,2) -- (5,2) -- (5,2.75) -- (1,2.75)--(1,2);

\draw [-stealth] (0,0)--(0.6,0);
\draw  (0.6,0)--(1,0);
\draw [-stealth] (0,1)--(0.6,1);
\draw  (0.6,1)--(1,1);
\draw [-stealth] (0,2)--(0.6,2);
\draw  (0.6,2)--(1,2);

\draw [-stealth] (0,0)--(0,0.6);
\draw  (0,0.6)--(0,1);
\draw [-stealth] (0,1)--(0,1.6);
\draw  (0,1.6)--(0,2);

\draw (0,2)--(0,2.3);
\fill (0,2.4) circle [radius=0.3mm];
\fill (0,2.5) circle [radius=0.3mm];
\fill (0,2.6) circle [radius=0.3mm];

\end{scope}

\fill (0,0) circle [radius=0.1];
\fill (1,0) circle [radius=0.1];
\fill (0,1) circle [radius=0.1];
\fill (1,1) circle [radius=0.1];
\fill (0,2) circle [radius=0.1];
\fill (1,2) circle [radius=0.1];

\node at (0.5,0) [below] {$a,c$};
\node at (0.5,1) [below] {$a,c$};
\node at (0.5,2) [below] {$a,c$};
\node at (0,0.5) [left] {$b,d$};
\node at (0,1.5) [left] {$b,d$};

\end{tikzpicture}

\caption{The graphs $\Gamma_M^{(3)}=\Cay(M,\{a,b\})$ and $\Gamma_N^{(3)}=\Cay(N,\{c,d\})$. Each copy of $\Gamma_M^{(2)}$ (resp. $\Gamma_N^{(2)}$) is indicated by its initial $a$ (resp. $c$) edge and a shaded rectangle.}
\label{fig3} \label{fig6}
\end{center}
\end{figure}

We begin with the description of $\Cay(M,\{a,b\})$.
The basic building block $\Gamma_M^{(1)}$ for this graph is shown in Figure \ref{fig1}.
In the view of the the relations $ab^na=aba$ a copy of $\Gamma_M^{(1)}$ will be found at every vertex
of $\Cay(M,\{a,b\})$ which receives an edge labelled $a$.
This process can be captured recursively as follows:
\begin{itemize}
\item
Start with a single edge labelled $a$.
\item
At each subsequent step, for every vertex $w$ introduced in the previous step with an in-edge labelled $a$ attach a new out-edge labelled $a$ and a copy of $\Gamma_M^{(1)}$ based at $w$.
\end{itemize}
The resulting graph $\Gamma_M^{(2)}$ is sketched in Figure \ref{fig2}.
This is the portion of $\Cay(M,\{a,b\})$ corresponding to the normal forms beginning with $a$.
All the normal forms are obtained by prepending these by an arbitrary power of $b$,
yielding the graph $\Gamma_M^{(3)}=\Cay(M,\{a,b\})$ shown in Figure \ref{fig3}.

Let us now turn our attention to $\Cay(N,\{c,d\})$.
The basic building block $\Gamma_N^{(1)}$, depicted in Figure \ref{fig4}, this time reflects the normal forms
$(d^3c)^q d^r$, and the fact that if they are pre- and post-multiplied by $c$ the resulting word `collapses' to $cdc$.
Since there will be a copy of this graph emanating from every vertex which has an in-edge labelled $c$, we can set up the following recursive process:
\begin{itemize}
\item
Start with a single edge labelled $c$.
\item
At each subsequent step, for every vertex $w$ introduced in the previous step with an in-edge labelled $c$ attach a new out-edge labelled $c$ and a a copy of $\Gamma_N^{(1)}$ based at $w$.
\end{itemize}
This results in the graph $\Gamma_N^{(2)}$ depicted in Figure \ref{fig5}.
It is important to notice that in this graph all the words of the form $u(d^3 c)^qd^r$
($u\in U_N$, $0\leq r\leq 3$) can be traced from the initial vertex, and that they all lead to distinct vertices.
Finally, by Lemma \ref{lemma2}, each of these normal forms can be preceded by an arbitrary power of $d$,
and so we obtain the graph $\Gamma_N^{(3)}=\Cay(N,\{c,d\})$ shown in Figure \ref{fig6}.

The above discussion and Figures \ref{fig1}, \ref{fig2}, \ref{fig3} 
show that when the edge labels are ignored the underlying graphs are isomorphic. We remark that the left Cayley graphs of $M$ and $N$ are not isomorphic, as is easily verified.

\section{Isomorphism}
\label{sec4}

We now provide a rigorous proof of our Main Theorem, by exhibiting an explicit isomorphism $f$ between
$\Cayu(M,\{a,b\}) $ and $\Cayu(N,\{c,d\})$.
Essentially, the mapping $f$ acts on a typical normal form $b^s ub^t$ ($u\in U_M$) by simple replacement
$a\mapsto c$, $b\mapsto d$ within the $b^su$ prefix, and transforms the final $b^t$ into
$(d^3c)^qd^r$ where $t=4q+r$, $0\leq r\leq 3$.
Since there are actually several types of normal forms, and since how a generator acts on a normal form depends on its type, the actual definition and the subsequent verification splits into cases.
For the sake of conciseness we have organised all the information in tabular form.

First of all, the normal forms for $M$ given by Lemma \ref{lemma1} are listed in Table \ref{table1}
alongside the permitted values of relevant parameters ($s$, $t$, $u$).
These values will subsequently be taken as read, without explicitly writing them down.
We take these normal forms for the vertices of $\Cayu(M,\{a,b\})$.
For each vertex there are two edges coming out of it, corresponding to multiplications by $a$ and $b$ respectively.
For instance, multiplying the normal form $b^sub^t$ of type \nf{NFM3} by $a$ and $b$ yields
\[
b^sub^t\cdot a\rightarrow^\ast b^suba,\ b^sub^t\cdot b=b^sub^{t+1},
\]
normal forms of types \nf{NFM2} and \nf{NFM3} respectively, resulting in two edges
\begin{equation}
\label{eq6}
(b^sub^t, b^suba), \ \ (b^sub^t, b^sub^{t+1}).
\end{equation}
The complete information on the edges of $\Cayu(M,\{a,b\})$ is given in Table \ref{table2}.
The only instance where the normal form differs from the simple concatenation with the generator is the row labelled NFM3a, where the rewriting rule $ab^ta\rightarrow aba$ is applied.

\begin{table}
  \begin{center}
  \begin{tabular}{ | l | l | l |}
    \hline
    Label & Normal Form & Parameters \\
    \hline
    \nf{NFM1} & $b^s$  &$s \geq 0$ \\
    \hline
    \nf{NFM2} & $b^su$ & $s \geq 0, u \in  U_M$ \\
    \hline
    \nf{NFM3} & $b^sub^t$ & $s \geq 0, u \in  U_M, t > 0$  \\
    \hline
  \end{tabular}  

  \caption{Normal forms for $M$} \label{table1}
\end{center}
\end{table}

\begin{table}
  \begin{center}
  \begin{tabular}{ | l | l | l |}
  \hline
  Label & Edge & Vertex Types \\
  \hline
  \nf{NFM1a}  & $(b^s, b^sa)$ & (\nf{NFM1},\nf{NFM2})\\
  \nf{NFM1b} & $(b^s, b^{s+1})$ & (\nf{NFM1},\nf{NFM1}) \\
  \hline
  \nf{NFM2a} & $(b^su,b^sua)$ & (\nf{NFM2},\nf{NFM2}) \\
  \nf{NFM2b} & $(b^su, b^sub)$ & (\nf{NFM2},\nf{NFM3})\\
  \hline
  \nf{NFM3a} & $(b^sub^t, b^suba)$ & (\nf{NFM3},\nf{NFM2})\\
  \nf{NFM3b}& $(b^sub^t, b^sub^{t+1})$ & (\nf{NFM3},\nf{NFM3}) \\
  \hline
  \end{tabular}   
  \caption{The edges in $\Cayu(M,\{a,b\})$} \label{table2}
\end{center}
\end{table}

The corresponding information for $\Cayu(N,\{c,d\})$ is given in Tables \ref{table3}, \ref{table4}.
It is worth noting that a further split in cases occurs for edges from vertices
$d^pv(d^3c)^qd^r$ of type \nf{NFN3}.
Indeed, bearing in mind that $v$ is non-empty and ends with $c$, depending on the value of $r$ we have:
\begin{align*}
&d^p v (d^3c)^q\cdot c = d^p v (d^3 c)^{q-1}d^3c^2\rightarrow d^pv(d^3c)^{q-1}dc\rightarrow^\ast d^pvdc
\\
&d^p v(d^3c)^qd\cdot c\rightarrow^\ast d^pvdc
\\
&d^pv(d^3c)^qd^2\cdot c\rightarrow d^pv(d^3c)^qdc\rightarrow^\ast d^pvdc
\\
&d^pv(d^3c)^qd^3\cdot c = d^pv(d^3c)^{q+1}.
\end{align*}
Thus, for $r=0,1,2$ we obtain normal forms of type \nf{NFN2}, while for $r=3$ it is \nf{NFN3}.
Similarly, for $r=0,1,2$ multiplication by $d$ increases $r$ by $1$, while for $r=3$ have
\[
d^pv(d^3c)^qd^3\cdot d=d^pv(d^3c)^qd^4\rightarrow d^pv(d^3c)^qdc\rightarrow^\ast d^pvdc.
\]

\begin{table} 
  \begin{center}
  \begin{tabular}{ | l | l | l | }
  \hline
  Label & Normal Form & Parameters \\
  \hline
  \nf{NFN1} & $d^p$ & $ p \geq 0$ \\
  \hline
  \nf{NFN2} & $d^pv$ & $p \geq 0, v \in  U_N$ \\
  \hline
  \nf{NFN3} & $d^pv(d^3c)^qd^r$ & $p,q \geq 0,\ v\in  U_N,\ 0 \leq r \leq 3,\ q+ r > 0$\\
  \hline
  \end{tabular}
  \caption{Normal forms for $N$} \label{table3}
\end{center}
\end{table}

\begin{table}
  \begin{center}
  \begin{tabular}{ | l | l | l |l|}
  \hline
  Label & Edge & Parameters & Vertex Types \\
  \hline
  \nf{NFN1c} & $(d^p, d^pc)$ && (\nf{NFN1},\nf{NFN2}) \\
  \nf{NFN1d}  & $(d^p,d^{p+1})$ && (\nf{NFN1},\nf{NFN1}) \\
  \hline
  \nf{NFN2c} & $(d^pv, d^pvc)$ && (\nf{NFN2},\nf{NFN2}) \\
  \nf{NFN2d} & $(d^pv, d^pvd)$ && (\nf{NFN2},\nf{NFN3})\\
  \hline
  \nf{NFN3c} & $(d^pv(d^3c)^qd^r, d^pvdc)$ & $0 \leq r \leq 2$ & (\nf{NFN3},\nf{NFN2})\\
     & $(d^pv(d^3c)^qd^3, d^pv(d^3c)^{q+1})$ && (\nf{NFN3},\nf{NFN3})\\
  \nf{NFN3d}  & $(d^pv(d^3c)^qd^r, d^pv(d^3c)^qd^{r+1})$ & $0 \leq r \leq 2$ & (\nf{NFN3},\nf{NFN3})\\
    & $(d^pv(d^3c)^qd^3, d^pvdc)$ && (\nf{NFN3},\nf{NFN2})\\ 
  \hline
  \end{tabular}  
  \caption{The edges in $\Cayu(N,\{c,d\})$} \label{table4}
\end{center}
\end{table}

We now define the mapping $f:\N_M\rightarrow \N_N$ which will turn out to be the desired isomorphism.
Note that there is a natural bijection $u\mapsto \overline{u}$ between $ U_M$ and $ U_N$
induced by the substitution $a\mapsto c$, $b\mapsto d$.
With this notation in mind, the mapping $f$ is defined in Table \ref{table5}.
It is clearly a bijection, and its inverse can be read off by reading the table from right to left.

\begin{table}
  \centering   
  \begin{tabular}{ | l | l | l | l | l | }
  \hline
  Type of $w$ & $w$ & Parameters & $f(w)$ & Type of $f(w)$\\
  \hline
  \nf{NFM1} & $b^s$  && $d^s$ & \nf{NFN1}\\
  \hline
  \nf{NFM2} & $b^su$ && $d^s\overline{u}$ & \nf{NFN2}\\
  \hline
  \nf{NFM3} & $b^sub^t$ & $t = 4q+ r$ & $d^s\overline{u}(d^3c)^qd^r$ & \nf{NFN3}\\
  \hline
  \end{tabular}  
  \caption{The definition of $f:\N_M\rightarrow \N_N$.} \label{table5}
\end{table}  

It remains to show that $f$ is a graph isomorphism, i.e. that both $f$ and $f^{-1}$ map edges to edges.
This is achieved by taking each edge type from Tables \ref{table2}, \ref{table4}, applying $f$ or $f^{-1}$ to its end-points using Table \ref{table5}, and verifying that the resulting pair of vertices also forms an edge.
For example, consider the edge
$(b^sub^t, b^suba)$ of type \nf{NFM3a} that we encountered in \eqref{eq6}. Under $f$ its endpoints are mapped to
$d^s \overline{u} (d^3c)^qd^r$, where $t=4q+r$, and $d^s\overline{u}dc$.
From Table \ref{table4} we see that for $r=0,1,2$, the pair $(d^s \overline{u} (d^3c)^qd^r, d^s\overline{u}dc)$ is an edge of type \nf{NFN3c}, while for $r=3$ it is still an edge but its type is \nf{NFN3d}.
The other edge that arose in \eqref{eq6} was $(b^sub^t, b^sub^{t+1})$ of type
\nf{NFM3b}. Its endpoints are mapped under $f$ to $d^s \overline{u} (d^3c)^qd^r$, where $t=4q+r$,
and $d^s\overline{u}(d^3 c)^{q_1}d^{r_1}$, where $t+1=4q_1+r_1$.
Noting that $q_1=q$ and $r_1=r+1$ for $r=0,1,2$, while $q_1=q+1$ and $r_1=0$ for $r=3$,
in the former case we obtain an edge of the type \nf{NFN3d}, and in the latter \nf{NFN3c}.
All the cases are presented in Tables \ref{table7}, \ref{table8}.

\begin{table}
  \begin{center}
  \begin{tabular}{ | l | l | l | l | l |}
  \hline
  Edge Type & $(w,v)$ & $(f(w), f(v))$ & Parameters & Edge Type \\
  \hline
  \nf{NFM1a} & $(b^s, b^sa)$  & $(d^s, d^sc)$ && \nf{NFN1c}\\
  \hline
  \nf{NFM1b} & $(b^s, b^{s+1})$& $(d^s, d^{s+1})$ && \nf{NFN1d}\\
  \hline
  \nf{NFM2a} & $(b^su, b^sua)$ & $(d^s\overline{u}, d^s\overline{u}c)$ && \nf{NFN2c} \\
  \hline
  \nf{NFM2b} & $(b^su, b^sub)$ & $(d^s\overline{u}, d^s\overline{u}d)$ && \nf{NFN2d} \\
  \hline
  \nf{NFM3a} & $(b^sub^t, b^suba)$ & $(d^s\overline{u}(d^3c)^qd^r,$ & $r=0,1,2$ & \nf{NFN3c}\\
   &&$\ \ \ \ \ \ \ \ \ d^s\overline{u}dc)$& $r=3$ & \nf{NFN3d}\\
  \hline
  \nf{NFM3b} & $(b^sub^t,  b^sub^{t+1})$ & $(d^s\overline{u}(d^3c)^q d^r,$ & $r=0,1,2$ & \nf{NFN3d}\\
&&$\ \ \  d^s\overline{u}(d^3c)^q d^{r+1})$  & & \\
  &  & $(d^s\overline{u}(d^3c)^q d^r, $ & $r = 3$ & \nf{NFN3c}\\
  &  & $\ \ \  d^s\overline{u}(d^3c)^{q+1})$ & & \\
  \hline
  \end{tabular}

  \caption{$f$ maps the edges of $\Cayu(M,\{a,b\})$ to edges of $\Cayu(N,\{c,d\})$} 
  \label{table7}
\end{center}
\end{table}

\begin{table}
\begin{center} 
  \begin{tabular}{| l | l | l | l | l |}
  \hline
  Edge Type & $(w, v)$ & $(f^{-1}(w), f^{-1}(v))$ & Parameters & Edge Type \\
  \hline
 \nf{NFN1c} & $(d^p, d^pc)$  & $(b^p, b^pa)$ && \nf{NFM1a} \\
  \hline
  \nf{NFN1d} & $(d^p, d^{p+1})$ & $(b^p, b^{p+1})$ && \nf{NFM1b} \\
  \hline
  \nf{NFN2c} & $(d^p\overline{u}, d^p\overline{u}c)$  & $(b^pu, b^pua)$ && \nf{NFM2a} \\
  \hline
  \nf{NFN2d} & $(d^p\overline{u}, d^p\overline{u}d)$ &  $(b^pu, b^pub)$ && \nf{NFM2b} \\
  \hline
  \nf{NFN3c} & $(d^p\overline{u}(d^3c)^qd^r,$    & $(b^pub^{4q+r} , b^puba)$ & $r=0,1,2$& \nf{NFM3a}\\
  & $\ \ \ \ \ \ \ \ \ d^p\overline{u}dc)$    &  & & \\
   & $(d^p\overline{u}(d^3c)^qd^3, $ & $(b^pub^{4q+3}, b^pub^{4q+4})$ && \nf{NFM3b}\\
& $\ \ d^p\overline{u}(d^3c)^{q+1})$ &  && \nf{NFM3b}\\
  \hline
  \nf{NFN3d} & $(d^p\overline{u}(d^3c)^qd^r,$   & $(b^pub^{4q+3}, b^pub^{4q+4})$ & $r=0,1, 2$ & \nf{NFM3b} \\
& $\ \ d^p\overline{u}(d^3c)^qd^{r+1})$   &  &  &  \\
   & $(d^p\overline{u}(d^3c)^qd^3, $ & $(b^pub^{4q+3}, b^puba)$ && \nf{NFM3a}\\
& $\ \ \ \ \ \ \ \ \ d^p\overline{u}dc$ & && \\
  \hline
  \end{tabular} 
  \caption{$f^{-1}$ maps the edges of $\Cayu(N,\{c,d\})$ to edges of $\Cayu(M,\{a,b\})$} 
\label{table8}
\end{center}
\end{table}

\medskip
\noindent
\textbf{Acknowledgement.}
The authors are grateful to an anonymous referee for their careful reading of the paper and suggestions which have improved the exposition.

\end{document}